\documentclass[12pt,titlepage]{article}
\usepackage{graphicx,amsmath,amssymb}

\numberwithin{equation}{section}
\renewcommand{\Pr}{\mathbf{P}}

\newtheorem{theo}{Theorem}[section]

\newcommand{\Sidak}{{\v{S}}id{\'a}k}
\newenvironment{proof}{\noindent {\bf Proof:}}{$\parallel$}
\newcommand{\MR}[1]{}

\begin{document}

\title{Classes of Multiple Decision Functions Strongly Controlling
FWER and FDR\thanks{The authors 
acknowledge support from
National Science Foundation (NSF) Grant DMS 0805809, 
National Institutes of Health (NIH) Grant RR17698, and
Environmental Protection Agency (EPA) Grant RD-83241902-0 to University of Arizona with
subaward number Y481344 to the University of South Carolina.}}

\author{Edsel A.\ Pe\~na\thanks{E.\ Pe\~na is Professor, Department of Statistics,
University of South Carolina, Columbia. E-Mail: pena@stat.sc.edu.}
\and Joshua D.\ Habiger\thanks{J.\ Habiger is a PhD student in the
Department of Statistics, University of South Carolina, Columbia. 
E-Mail: habigerj@mailbox.sc.edu.} 
\and Wensong Wu\thanks{W.\ Wu is a PhD student in the
Department of Statistics, University of South Carolina, Columbia.
E-Mail: wu26@mailbox.sc.edu.} \\ \medskip \and
Department of Statistics \\ 
University of South Carolina \\
Columbia, SC 29208 USA}

\maketitle

\begin{abstract}
This paper provides two general classes of multiple decision functions where each member of
the first class strongly controls the family-wise error rate (FWER), while each member
of the second class strongly controls the false discovery rate (FDR).
These classes offer the possibility that an
optimal multiple decision function with respect to a
pre-specified criterion, such as the
missed discovery rate (MDR), could be found within these classes.
Such multiple decision functions can be utilized in multiple testing,
specifically, but not limited to, the analysis of
high-dimensional microarray data sets.
\end{abstract}

\medskip

\noindent
{\bf Keywords and Phrases:} 
{false discovery rate};
{family wise error rate};
{missed discovery rate};
{multiple decision problem};
{multiple testing};
{strong control.}

\section{Introduction}
\label{section-Introduction}

Consider the situation which arises in the analysis
of high-dimensional data, epitomized by a microarray data, where $M$
pairs of null and alternative hypotheses, $(H_{m0},H_{m1}), m=1,2,\ldots,M$,
are simultaneously tested; see, for instance, \cite{Efr08,EfrAAS08} 
for concrete examples of such situations. Two commonly-used Type I error
rates for this multiple testing problem are the family-wise error rate (FWER),
which is the probability of at least one false discovery, where discovery
means rejecting (accepting) a null (an alternative) hypothesis, and the false
discovery rate (FDR), which is the expectation of the ratio of the number
of false discoveries over the number of discoveries. The usual testing
paradigm employed in these situations is to decide on the collection of
statistical tests for
the $M$ pairs of hypotheses, e.g., a $t$-test for each pair, obtain the
$p$-value for each test, and then use the resulting $M$ $p$-values 
in the FWER-controlling 
sequential \Sidak\ procedure, provided an independence condition
is satisfied, or the FDR-controlling procedure in \cite{BenHoc95}. In this
conventional approach, there appears to be no leeway in the choice of the multiple
testing procedure the moment the individual test
procedures have been chosen. 

However, we pose the following question.
If we are given the $M$ test procedures for each of the $M$ pairs of hypotheses,
could we obtain classes of multiple testing procedures whose elements
either control the FWER or the FDR? If the answer to this
question is in the affirmative, then
we may be able to find {\em a}  multiple testing procedure within these classes
which is optimal with respect to some chosen Type II error rate. 
And, we may then be able to
choose the starting collection of test functions that will provide the best
multiple testing procedure. 

This paper is in this spirit. We will demonstrate
that, under certain conditions, when given a collection of test functions for
the $M$ pairs of hypotheses, that we can generate classes of multiple
testing procedures controlling the FWER or the FDR. The results have important
implications in the search for optimal multiple testing procedures that
control either of these Type I error rates as we will see later.
We shall investigate these issues in a general, but not surprisingly, 
more abstract framework. The main results in this paper were motivated
by those in \cite{PenHabWu10,PenHabWu10Supp} which did not deal with
classes of multiple testing procedures, but instead focussed in developing
improved FWER and FDR-controlling procedures from the Neyman-Pearson most
powerful tests for each of the $M$ pairs of hypotheses.

\section{Mathematical Setting}
\label{section-Setting}

Let $(\mathcal{X},\mathcal{F},\mathcal{P})$ be a statistical model,
so $(\mathcal{X},\mathcal{F})$ is a measurable space and
$\mathcal{P}$ is a collection of probability measures on
$(\mathcal{X},\mathcal{F})$. Though not needed
in the abstract development, for concreteness we may adopt the 
usual interpretation that $\mathcal{X}$ is the space of possible realizations of an
observable random entity $X$ from an experiment or a study.
In decision problems with action space
$\mathfrak{A} = \{0,1\}$, such as in hypothesis testing,
a nonrandomized decision function is a $\delta: (\mathcal{X},\mathcal{F})
\rightarrow (\mathfrak{A},\sigma(\mathfrak{A}))$.
In the hypothesis testing setting, given $X = x \in \mathcal{X}$,
a decision $\delta(x) = 0$ corresponds to
deciding in favor of a null hypothesis ($H_0$), 
whereas a decision of $\delta(x) = 1$, a {\em discovery},
corresponds to rejecting $H_0$ in favor of an alternative hypothesis ($H_1$).

It suffices to restrict ourselves to 
nonrandomized decision functions since, through the use of an auxiliary randomizer which
is usually a standard uniform variable $U$ that is independent of 
$X$, we can always convert a randomized decision function $\delta^*:
(\mathcal{X},\mathcal{F}) \rightarrow ([0,1],\sigma[0,1])$
into a nonrandomized decision function 
$\delta: (\mathcal{X} \times [0,1],\mathcal{F} \otimes \sigma[0,1])
\rightarrow (\mathfrak{A},\sigma(\mathfrak{A}))$
via $\delta(x,u) = I\{u \le \delta^*(x)\}$ with $I\{\cdot\}$ the indicator
function.
 
Thus, in our general formulation, the sample space $\mathcal{X}$ may actually
represent a product space between a data space and $[0,1]$. This framework
is appropriate, for instance, when dealing with discrete data or when using
nonparametric decision functions. For more discussions on this matter, see
\cite{PenHabWu10,HabPen10}. 

Decision or test functions typically depend on a size parameter
$\alpha \in [0,1]$. For example, when testing the null hypothesis
$H_0: \mu = 0$ versus the alternative hypothesis $H_1: \mu \ne 0$
based on a random observable $X \sim N(\mu,1)$, the size-$\alpha$ test 
$\delta: \mathcal{X} \equiv \Re \rightarrow \{0,1\}$ has
%
$\delta(x;\alpha) = I\{|x| > \Phi^{-1}(1-\alpha/2)\},$
%
where $\Phi^{-1}(\cdot)$ is the quantile function of the standard normal
distribution. Henceforth, in order to simplify our notation,
we shall adopt a functional notation where
$\delta(\alpha)$ represents the statistic defined on $\mathcal{X}$ 
according to $x \mapsto \delta(x;\alpha)$.
Now, when viewed as a process in $\alpha$, we then obtain the notion
of a (nonrandomized) decision process introduced in \cite{PenHabWu10}, which is
a stochastic process $\Delta = \{\delta(\alpha): \ \alpha \in [0,1]\}$
where, $\forall \alpha \in [0,1]$, $\delta(\alpha)$ is a decision function,
and such that the following conditions are satisfied.
\begin{itemize}
\item[(D1)] $\delta(0) = 0$ and $\delta(1) = 1$ a.e.-$\mathcal{P}$.
\item[(D2)] The sample paths $\alpha \mapsto \delta(\alpha)$
are, a.e.-$\mathcal{P}$, $\{0,1\}$-valued step-functions which are nondecreasing and
right-continuous.
\end{itemize}

Let $\mathcal{M}$ be a finite set with $|\mathcal{M}| = M$. An
$\mathcal{M}$-indexed multiple decision problem is one whose
action space is $\mathfrak{A}^M$. In the context of a multiple
hypotheses testing problem, for each $m \in \mathcal{M}$, there
is a pair of hypotheses $H_{m0}$ and $H_{m1}$. Of interest is to
simultaneously decide between $H_{m0}$ and $H_{m1}$ for each
$m \in \mathcal{M}$.
A multiple decision function (MDF) for such a problem is a
$\mathbf{\delta} = (\delta_m: \ m \in \mathcal{M})$ where $\delta_m$
is a decision function, so that $\delta: (\mathcal{X},\mathcal{F})
\rightarrow (\mathfrak{A}^M,\sigma(\mathfrak{A}^M))$.
A multiple decision process (MDP) is a
$\mathbf{\Delta} = (\Delta_m:\ m \in \mathcal{M})$
where $\Delta_m = \{\delta_m(\alpha):\ \alpha \in [0,1]\}$ is a decision process.

For each $\Pr \in \mathcal{P}$, let there be subsets
$\mathcal{M}_0(\Pr)$ and $\mathcal{M}_1(\Pr)$ of $\mathcal{M}$
such that
\begin{displaymath}
\mathcal{M} = \mathcal{M}_0(\Pr) \cup \mathcal{M}_1(\Pr) \quad \mbox{and} \quad
\mathcal{M}_0(\Pr) \cap \mathcal{M}_1(\Pr) = \emptyset.
\end{displaymath}
We shall assume that the following condition holds.
\begin{itemize}
\item[(D3)] Under $\Pr$, the subcollections $\{\Delta_m:\ m \in \mathcal{M}_0(\Pr)\}$ and
$\{\Delta_m:\ m \in \mathcal{M}_1(\Pr)\}$ are independent of each other,
and the elements of $\{\Delta_m:\ m \in \mathcal{M}_0(\Pr)\}$ are independent.
\end{itemize}
In the multiple hypotheses testing
situation, $H_{m0}$ is true under $\Pr$ if and only if $m \in \mathcal{M}_0(\Pr)$.
Observe that the elements
of $\{\Delta_m:\ m \in \mathcal{M}_1(\Pr)\}$
need not be independent of each other, under $\Pr$.
We shall also assume that
\begin{itemize}
\item[(D4)] \label{size}
With $E_\Pr(\cdot)$ denoting the expectation operator under $\Pr$,
then $\forall \Pr \in \mathcal{P},
\forall m \in \mathcal{M}_0(\Pr),
\forall \alpha \in [0,1]$, we have $E_\Pr \left\{ \delta_m(\alpha) \right\} = \alpha.$
\end{itemize}
The collection of all $\mathcal{M}$-indexed multiple decision processes
satisfying conditions (D1)--(D4) will be denoted by $\mathfrak{D}$. We remark
that the requirement of equality in (D4) given by
$E_\Pr \left\{ \delta_m(\alpha) \right\} = \alpha$
will usually be fulfilled in many situations since an auxiliary randomizer is
incorporated in our framework, though there may still
be situations when dealing with non-regular families of distributions 
where this condition may not be satisfied. The latter will manifest itself when
the decision functions already have power equal to one but without yet requiring
their sizes to equal one.

Let $\mathbf{A} = (A_m:\ m \in \mathcal{M})$ be an
$\mathcal{M}$-indexed collection of measurable
functions with $A_m: ([0,1],\sigma[0,1]) \rightarrow ([0,1],\sigma[0,1])$. We
assume that, for each $m \in \mathcal{M}$, the following conditions are
satisfied:
\begin{itemize}
\item[(A1)]
$A_m(0) = 0$ and $A_m(1) = 1$.
\item[(A2)]
The mapping $\alpha \mapsto A_m(\alpha)$ is continuous and strictly increasing.
\item[(A3)]
$\forall \alpha \in [0,1],\ \prod_{m\in\mathcal{M}} [1 - A_m(\alpha)] \ge
1 - \alpha$.
\item[(A4)] \label{size condition}
$\forall \alpha \in [0,1], \forall \Pr \in \mathcal{P}:
|\mathcal{M}_0(\Pr)| \max_{m \in \mathcal{M}_0(\Pr)} A_m(\alpha) \le
\sum_{m \in \mathcal{M}} A_m(\alpha).$
%
%
\end{itemize}
Such an $\mathbf{A}$ will be called a multiple decision size function.
The collection of all $\mathcal{M}$-indexed multiple decision size functions
will be denoted by $\mathfrak{S}$. A particular element of
$\mathfrak{S}$ is the Sidak multiple decision size function (cf., \cite{Sid67})
$\mathbf{A}^S = (A_m^S:\ m \in \mathcal{M})$ with
\begin{equation}
\label{Sidak size}
A_m^S(\alpha) = 1 - (1 - \alpha)^{{1}/{M}},\ \alpha \in [0,1], m \in \mathcal{M}.
\end{equation}
Another particular element of $\mathfrak{S}$ is the Bonferroni size function
$\mathbf{A}^B = (A_m^B:\ m \in \mathcal{M})$ with
\begin{equation}
\label{Bonferroni size}
A_m^B(\alpha) = \alpha/M, \ \alpha \in [0,1], m \in \mathcal{M}.
\end{equation}

Before proceeding we also recall the notion of generalized $P$-value statistics;
see \cite{PenHabWu10}.
Given a $\mathbf{\Delta} \in \mathfrak{D}$ and an $\mathbf{A} \in \mathfrak{S}$,
we define for $m \in \mathcal{M}$ the random variable
\begin{equation}
\alpha_m \equiv \alpha_m(\mathbf{\Delta},\mathbf{A}) =
\inf\left\{
\alpha \in [0,1]: \delta_m(A_m(\alpha)) = 1
\right\}.
\label{generalized p-value}
\end{equation}
The collection $(\alpha_m(\mathbf{\Delta},\mathbf{A}): m \in \mathcal{M})$
is called the vector of generalized $P$-value statistics associated with the pair
$(\mathbf{\Delta},\mathbf{A})$. Observe that the usual $P$-value statistic associated
with $\delta_m$ is $P_m = A_m(\alpha_m)$, hence the use of the adjective
{\em generalized} for the $\alpha_m$s. We shall assume without much loss of
generality that these generalized $P$-values are a.e.\ $[\mathcal{P}]$ distinct.

\section{Main Theorems and Classes of MDFs}
\label{section-Main Theorems}

We shall present in this section the two main results that will enable the
construction of the classes of multiple decision functions controlling
FWER and FDR.

Given a $\mathbf{\Delta} = \{\Delta_m:\ m \in \mathcal{M}\} \in \mathfrak{D}$,
an $\mathbf{A} = \{A_m:\ m \in \mathcal{M}\} \in \mathfrak{S}$, a
$\Pr \in \mathcal{P}$, and an $\alpha \in [0,1]$, define the stochastic
processes $\mathbf{S}_0 = \{S_0(\alpha): \alpha \in [0,1]\}$,
$\mathbf{S} = \{S(\alpha): \alpha \in [0,1]\}$, and
$\mathbf{F} = \{F(\alpha): \alpha \in [0,1]\}$, where
\begin{eqnarray}
S_0(\alpha) & \equiv &
S_0(\alpha;\mathbf{\Delta},\mathbf{A},\Pr) =
\sum_{m \in \mathcal{M}_0(\Pr)} \delta_m(A_m(\alpha)); \label{S0} \\
S(\alpha) & \equiv &
S(\alpha;\mathbf{\Delta},\mathbf{A}) =
\sum_{m \in \mathcal{M}} \delta_m(A_m(\alpha));  \label{S} \\
F(\alpha) & \equiv &
F(\alpha;\mathbf{\Delta},\mathbf{A},\Pr) =
\frac{S_0(\alpha)}{S(\alpha)} I\{S(\alpha) > 0\}, \label{F}
\end{eqnarray}
with the convention that ${0}/{0} = 0$.
These quantities have the following interpretations. 
Given an $\alpha \in [0,1]$, for each $m \in \mathcal{M}$,
the decision function whose size is $A_m(\alpha)$ is chosen from
$\Delta_m$, and the MDF $\delta(\alpha) \equiv (\delta_m[A_m(\alpha)]: m \in
\mathcal{M})$ will be employed in the decision-making. For this MDF $\delta(\alpha)$,
then $S_0(\alpha)$ is the number of false discoveries, $S(\alpha)$ is
the number of discoveries, and $F(\alpha)$ is the proportion of false
discoveries among all discoveries.
Observe, however, that since $\Pr$ is unknown, both $\mathbf{S}_0$ and
$\mathbf{F}$ are unobservable, whereas $\mathbf{S}$ is observable.

For $q \in [0,1]$, let us also define the random variables
\begin{eqnarray}
\lefteqn{ \alpha^\dagger(q)  \equiv  \alpha^\dagger(q;\mathbf{\Delta},\mathbf{A}) 
\nonumber } \\ & = &
\inf
\left\{
\alpha \in [0,1]:
\prod_{m \in \mathcal{M}}
[1 - A_m(\alpha)]^{1-\delta_m(A_m(\alpha)-)} < 1 - q
\right\}; \label{alpha-dag}
\end{eqnarray}
and
\begin{eqnarray}
\lefteqn{ \alpha^*(q) \equiv \alpha^*(q;\mathbf{\Delta},\mathbf{A}) \nonumber } \\ & = &
\sup\left\{
\alpha \in [0,1]:\
\sum_{m \in \mathcal{M}} A_m(\alpha) \le q S(\alpha;\mathbf{\Delta},\mathbf{A})
\right\}. \label{alpha-star}
\end{eqnarray}
In essence, $\alpha^\dagger(q)$ is a {\em first} crossing-time random variable, 
whereas $\alpha^*(q)$ is a {\em last} crossing-time random variable.
The forms of these two random variables were motivated and
justified in Sections 6 and 7 in \cite{PenHabWu10} for
a specific multiple decision size function, but the justifications
in that paper carry over to the more general setting considered
here.

The two main results of this paper are contained in
Theorem \ref{theo-FWER control} and Theorem \ref{theo-FDR control}.
We present the statements of these theorems, but defer their proofs
to Section \ref{section-Proof} after some discussions about their
implications and potential usefulness.

\begin{theo}
\label{theo-FWER control}
Under conditions (D1)--(D4) for $\mathfrak{D}$ and (A1)--(A3) for $\mathfrak{S}$, 
we have that
$\forall \Pr \in \mathcal{P},
\forall \mathbf{\Delta} \in \mathfrak{D},
\forall \mathbf{A} \in \mathfrak{S},
\forall q \in [0,1]$,
\begin{displaymath}
E_\Pr\left\{I\{S_0(\alpha^\dagger(q;\mathbf{\Delta},\mathbf{A});
\mathbf{\Delta},\mathbf{A},\Pr) \ge 1\}\right\} \le q.
\end{displaymath}
\end{theo}

Observe that $E_\Pr\left\{I\{S_0(\alpha^\dagger(q;\mathbf{\Delta},\mathbf{A});
\mathbf{\Delta},\mathbf{A},\Pr) \ge 1\}\right\}$ is the FWER since it is
the probability of
committing at least one false discovery when the true underlying probability measure
is $\Pr$. Thus, Theorem \ref{theo-FWER control} shows that for any $q \in [0,1]$,
any multiple decision process $\mathbf{\Delta} \in \mathfrak{D}$,
and any multiple decision size function $\mathbf{A} \in \mathfrak{S}$, 
the MDF defined via
\begin{equation}
\label{FWER-controlling MDF}
\mathbf{\delta}^\dagger(q) \equiv
\mathbf{\delta}^\dagger(q;\mathbf{\Delta},\mathbf{A}) =
\left(
\delta_m[A_m(\alpha^\dagger(q;\mathbf{\Delta},\mathbf{A}))]:\ m \in \mathcal{M}
\right),
\end{equation}
{\em strongly} controls the FWER at $q$.

\begin{theo}
\label{theo-FDR control}
Under conditions (D1)--(D4) for $\mathfrak{D}$ and (A1)--(A4) for $\mathfrak{S}$,
we have that $\forall \Pr \in \mathcal{P},
\forall \mathbf{\Delta} \in \mathfrak{D},
\forall \mathbf{A} \in \mathfrak{S},
\forall q \in [0,1]$,
\begin{displaymath}
E_\Pr\left\{F(\alpha^*(q;\mathbf{\Delta},\mathbf{A});
\mathbf{\Delta},\mathbf{A},\Pr)\right\} \le q.
\end{displaymath}
\end{theo}

Note that $E_\Pr\left\{F(\alpha^*(q;\mathbf{\Delta},\mathbf{A});
\mathbf{\Delta},\mathbf{A},\Pr)\right\}$
is the FDR as introduced in the seminal paper of \cite{BenHoc95}.
The implication of Theorem \ref{theo-FDR control} is that if, for each $q \in [0,1]$, and for
any multiple decision process $\mathbf{\Delta} \in \mathfrak{D}$ and 
multiple decision size function $\mathbf{A} \in \mathfrak{S}$, we define
the MDF
\begin{equation}
\label{FDR-control MDF}
\mathbf{\delta}^*(q) \equiv
\mathbf{\delta}^*(q;\mathbf{\Delta},\mathbf{A}) =
\left(
\delta_m[A_m(\alpha^*(q;\mathbf{\Delta},\mathbf{A}))]:\ m \in \mathcal{M}
\right),
\end{equation}
then $\mathbf{\delta}^*(q)$ is an MDF that
controls the FDR at $q$.

The importance of the preceding results is that
each multiple decision process $\mathbf{\Delta} \in \mathfrak{D}$
may have an associated multiple decision size process
$\mathbf{A} \equiv \mathbf{A}(\mathbf{\Delta}) \in
\mathfrak{S}$ such that the resulting multiple decision functions
$\mathbf{\delta}^\dagger(q)$ or
$\mathbf{\delta}^*(q)$ possess some optimality property, for example,
with respect to the missed discovery rate. To define this rate, let
\begin{equation}
\label{M}
M(\alpha)  \equiv
M(\alpha;\mathbf{\Delta},\mathbf{A},\Pr) =
\frac{\sum_{m \in \mathcal{M}_1(\Pr)} (1-\delta_m(A_m(\alpha)))}
{|\mathcal{M}_1(\Pr)|} I\{|\mathcal{M}_1(\Pr)| > 0\}.
\end{equation}
The quantity $M(\alpha)$ has the interpretation of being the proportion
of missed discoveries relative to the number of correct alternative
hypotheses. Then, for instance,
the missed discovery rate (MDR) of the MDF
in (\ref{FDR-control MDF}) is
\begin{displaymath}
E_\Pr\left\{M(\alpha^*(q); \mathbf{\Delta},\mathbf{A},\Pr)\right\}.
\end{displaymath}
For the given $\mathbf{\Delta}$, with proper choice of $\mathbf{A}$,
we may be able to find an MDF that strongly
controls the FWER or the FDR,
while at the same time possessing an optimal property with respect to another criterion,
such as having a small, possibly maximally over $\mathcal{P}$, MDR. 
This idea was implemented in a more restricted setting in
\cite{PenHabWu10,PenHabWu10Supp} when each of the pairs of
hypotheses contained
simple null and simple alternative hypotheses.

We note that previous works usually
focussed in developing {\em a} particular MDF and then
verifying that it controls the FWER or the FDR, such as, for example,
in \cite{BenHoc95}; more comprehensively, see \cite{DudLaa08}. 
It is our hope that by providing a class of
MDFs where each member strongly controls the FWER, given by
\begin{equation}
\label{class of FWER-control}
\mathfrak{C}^\dagger =
\left\{
\mathbf{\delta}^\dagger(q;\mathbf{\Delta},\mathbf{A}):\
\mathbf{\Delta} \in \mathfrak{D},
\mathbf{A} \in \mathfrak{S}
\right\};
\end{equation}
or a class of MDFs where each member controls the FDR, given by
\begin{equation}
\label{class of FDR-control}
\mathfrak{C}^* =
\left\{
\mathbf{\delta}^*(q;\mathbf{\Delta},\mathbf{A}):\
\mathbf{\Delta} \in \mathfrak{D},
\mathbf{A} \in \mathfrak{S}
\right\},
\end{equation}
then we will acquire the possibility of selecting from these classes 
MDFs possessing other desirable properties with respect to Type II error
rates. More discussion of this issue will be provided in
Section \ref{section-Optimal MDFs}.


\section{Proofs of the Main Theorems}
\label{section-Proof}

The proofs of the two theorems are analogous to those of
Theorem 6.1 and Theorem 7.1 in \cite{PenHabWu10} which can be found in the
supplemental article \cite{PenHabWu10Supp}. Note that those proofs were for
special forms of the multiple decision process and
multiple decision size function, whereas in the current
paper we are dealing with an arbitrary element $\mathbf{\Delta} \in \mathfrak{D}$
and an arbitrary element $\mathbf{A} \in \mathfrak{S}$. In the proofs
below, we assume that $\mathbf{\Delta} \in \mathfrak{D}$ 
and $\mathbf{A} \in \mathfrak{S}$ have been chosen and are fixed. Also, $q \in [0,1]$ and
$\Pr \in \mathcal{P}$ denotes the unknown underlying probability measure.
The dependence on $(\mathbf{\Delta}, \mathbf{A}, \Pr)$ of some of the relevant processes
and quantities below will not be explicitly written for brevity, unless needed
for clarity.

\subsection{Of Theorem \ref{theo-FWER control}}
\label{subsection-Proof of First Theorem}

\begin{proof}
We start by defining the stochastic process 
$\mathbf{H}_1 = \{H_1(\alpha): \alpha \in [0,1]\}$ via
\begin{equation}
\label{H1}
H_1(\alpha) \equiv H_1(\alpha;\mathbf{\Delta},\mathbf{A}) =
\prod_{m \in \mathcal{M}} [1 - A_m(\alpha)]^{1 - \delta_m(A_m(\alpha)-)}.
\end{equation}
The sample paths of this process are, a.e.\ $[\Pr]$,
left-continuous with right-hand limits ({\em caglad}),
are piecewise nonincreasing, and with 
$$1 - \alpha \le H_1(\alpha-) = H_1(\alpha) \le H_1(\alpha+)$$
for every $\alpha \in (0,1)$, where the first inequality
is due to property (A3). In fact, by virtue of property (A1) and
property (D1), note that
\begin{displaymath}
\lim_{\alpha \downarrow 0} H_1(\alpha) = 1 \quad
\mbox{and} \quad
\lim_{\alpha \uparrow 1} H_1(\alpha) = 1.
\end{displaymath}
Now, in terms of $\mathbf{H}_1$, we have that
\begin{displaymath}
\alpha^\dagger(q) = \inf\left\{\alpha \in [0,1]: H_1(\alpha) < 1 - q\right\}.
\end{displaymath}
Since, as pointed out above, we have $1 - \alpha \le H_1(\alpha)$, then
by its definition, we must have $\alpha^\dagger(q) \ge q$. This implies that
\begin{equation}
\label{H1 at alpha dag}
H_1(\alpha^\dagger(q)) \ge 1 - q.
\end{equation}

For the quantity of main interest in the theorem, we now have 
\begin{eqnarray*}
\lefteqn{
E_\Pr \left[
I\left\{ S_0(\alpha^\dagger(q)) \ge 1 \right\} \right] } \\ & = &
\Pr\left\{ S_0(\alpha^\dagger(q)) \ge 1 \right\} \\
& = &  1 - \Pr\left\{ S_0(\alpha^\dagger(q)) = 0 \right\} \\
& = & 1 - \Pr\left\{
\bigcap_{m \in \mathcal{M}_0(\Pr)}
\left[
\delta_m(A_m(\alpha^\dagger(q))) = 0
\right]\right\}.
\end{eqnarray*}
The last probability cannot, however, be written as a product of
probabilities since the $\delta_m(A_m(\alpha^\dagger(q)))$ for
$m \in \mathcal{M}_0(\Pr)$ need not be independent owing to the
dependence on $\alpha^\dagger(q)$ which is determined by all the
$(\Delta_m, m \in \mathcal{M})$.
On the other hand, we do have the set equality
\begin{equation}
\bigcap_{m \in \mathcal{M}_0(\Pr)}
\left[
\delta_m(A_m(\alpha^\dagger(q))) = 0
\right] =
\left\{
\alpha^\dagger(q) < \min_{m \in \mathcal{M}_0(\Pr)} \alpha_m
\right\},
\label{eq1}
\end{equation}
where the $\alpha_m$s are the generalized $p$-value statistics defined
in (\ref{generalized p-value}).

Next, define the stochastic process $\mathbf{H}_2 = \{H_2(\alpha): \alpha \in [0,1]\}$ via
\begin{eqnarray*}
\lefteqn{ H_2(\alpha)  \equiv  H_2(\alpha;\mathbf{\Delta},\mathbf{A},\Pr) } \\ & = &
\left(
\prod_{m \in \mathcal{M}_0(\Pr)}
[1 - A_m(\alpha)]
\right)
\left(
\prod_{m \in \mathcal{M}_1(\Pr)}
[1 - A_m(\alpha)]^{1 - \delta_m(A_m(\alpha)-)}
\right).
\end{eqnarray*}
Analogously to the $\mathbf{H}_1$ process, this has caglad sample paths.
Let us then define the quantity
\begin{eqnarray*}
\alpha^{\#}(q)  \equiv
\alpha^{\#}(q;\mathbf{\Delta},\mathbf{A},\Pr)  =
\inf\left\{
\alpha \in [0,1]:
H_2(\alpha)
< 1 - q
\right\}.
\end{eqnarray*}
Note that this is not a random variable since this depends on the unknown
probability measure $\Pr$, in contrast to $\alpha^\dagger(q)$. Furthermore,
also note that
\begin{equation}
\label{H2 at alpha pound}
H_2(\alpha^{\#}(q)) \ge 1-q.
\end{equation}
From their definitions, $H_1(\alpha) \ge H_2(\alpha)$, so that
$H_1(\alpha) < 1 - q$ implies $H_2(\alpha) < 1 - q$. Consequently,
\begin{equation}
\label{relationship between alpha-dag and alpha-pound}
\alpha^\dagger(q) \ge \alpha^\#(q).
\end{equation}
Now, the importance of the quantity $\alpha^{\#}(q)$ arises because of the
crucial set equality
\begin{equation}
\left\{
\alpha^\dagger(q) < \min_{m \in \mathcal{M}_0(\Pr)} \alpha_m
\right\}
=
\left\{
\alpha^{\#}(q) < \min_{m \in \mathcal{M}_0(\Pr)} \alpha_m
\right\}.
\label{eq2}
\end{equation}
To see this equality, first observed that the inclusion $\subseteq$
immediately follows from (\ref{relationship between alpha-dag and alpha-pound}).
To prove the reverse inclusion, since 
$$\{\alpha^\#(q) < \min_{m \in \mathcal{M}_0(\Pr)} \alpha_m\}$$
implies that, for some $\alpha_0 < \min_{m \in \mathcal{M}_0(\Pr)} \alpha_m$,
we have
$H_2(\alpha_0) < 1 - q.$
But for such an $\alpha_0$, we will have $\delta_m(A_m(\alpha_0)-) = 0$
for all $m \in \mathcal{M}_0(\Pr)$, so that
$$\alpha_0 \in \{\alpha \in [0,1]: H_1(\alpha) < 1 - q\}.$$
Consequently,
$$\alpha^\dagger(q) = \inf\{\alpha \in [0,1]: H_1(\alpha) < 1 - q\}
\le \alpha_0 < \min_{m \in \mathcal{M}_0(\Pr)} \alpha_m.$$
The reverse inclusion $\supseteq$ thus follows, completing the proof
of (\ref{eq2}).

By (\ref{eq1}), (\ref{eq2}), and the iterated expectation rule, it now
follows that
\begin{eqnarray*}
\lefteqn{ \Pr\left\{
\bigcap_{m \in \mathcal{M}_0(\Pr)}
\left[
\delta_m(A_m(\alpha^\dagger(q))) = 0
\right]
\right\} } \\
& = & \Pr\{\alpha^\dagger(q) < \min_{m \in \mathcal{M}_0(\Pr)} \alpha_m\}
 \\
& = &
\Pr\left\{
\alpha^{\#}(q) < \min_{m \in \mathcal{M}_0(\Pr)} \alpha_m
\right\} \\
& = &
E_\Pr
\left[
\Pr\left\{\alpha^{\#}(q) < \min_{m \in \mathcal{M}_0(\Pr)} \alpha_m
\ \left|\  \alpha^{\#}(q) \right. \right\}
\right].
\end{eqnarray*}
Since $\alpha^{\#}(q)$ is measurable with respect to the sub-$\sigma$-field
$\sigma(\delta_m:
m \in \mathcal{M}_1(\Pr))$, whereas $\min_{m \in \mathcal{M}_0(\Pr)} \alpha_m$
is measurable with respect to the sub-$\sigma$-field
$\sigma(\delta_m: m \in \mathcal{M}_0(\Pr))$, then
by condition (D3), $\alpha^{\#}(q)$ and $\min_{m \in \mathcal{M}_0(\Pr)} \alpha_m$
are independent. Furthermore, by condition (D3), we obtain
\begin{eqnarray*}
\Pr\left\{
\min_{m \in \mathcal{M}_0(\Pr)} \alpha_m > w
\right\}
& = &
\Pr\left\{
\bigcap_{m \in \mathcal{M}_0(\Pr)}
\left[
\delta_m(A_m(w)) = 0
\right]
\right\}
\\
& = &
\prod_{m \in \mathcal{M}_0(\Pr)}
\Pr\{\delta_m(A_m(w)) = 0\} \\
& = &
\prod_{m \in \mathcal{M}_0(\Pr)}
\left[
1 - A_m(w)
\right],
\end{eqnarray*}
with the last equality a consequence of condition (D4).  
Therefore,
\begin{eqnarray*}
\lefteqn{
\Pr\left\{
\bigcap_{m \in \mathcal{M}_0(\Pr)}
\left[
\delta_m(A_m(\alpha^\dagger(q))) = 0
\right]\right\} } \\ &  = &
E_\Pr
\left\{
\prod_{m \in \mathcal{M}_0(\Pr)}
\left[
1 - A_m(\alpha^{\#}(q))
\right]\right\}  \\
& \ge &
E_\Pr
\left\{
\left(
\prod_{m \in \mathcal{M}_0(\Pr)}
\left[
1 - A_m(\alpha^{\#}(q))
\right]
\right) \times \right. \\ && \left.
\left(
\prod_{m \in \mathcal{M}_1(\Pr)}
\left[
1 - A_m(\alpha^{\#}(q))
\right]^{1-\delta_m(A_m(\alpha^{\#}(q))-)}
\right)
\right\} \\
& = & E_\Pr\{H_2(\alpha^{\#}(q)\} \\
& \ge & E_\Pr (1 - q) \\ & = & 1-q
\end{eqnarray*}
with the last inequality following from (\ref{H2 at alpha pound}).
Thus, finally, we have
\begin{eqnarray*}
E_\Pr
\left[
I\left\{
S_0(\alpha^\dagger(q)) \ge 1
\right\}
\right]
& = & 1 - \Pr\left\{
\bigcap_{m \in \mathcal{M}_0(\Pr)}
\left[
\delta_m(A_m(\alpha^\dagger(q))) = 0
\right]\right\} \\
& \le & 1 - (1 - q) \\ &  = & q.
\end{eqnarray*}
This completes the proof of Theorem \ref{theo-FWER control}.
\end{proof}

We remark that condition (D4) can be weakened to just having
\begin{equation}
\label{weaker D4}
\forall m \in \mathcal{M}_0(\Pr), \forall \alpha \in [0,1]:\ E_\Pr\{\delta_m(\alpha)\} \le \alpha
\end{equation}
to still get the desired strong FWER control. This is so since in the
portion of the proof where we have
\begin{displaymath}
\prod_{m \in \mathcal{M}_0(\Pr)}
\Pr\{\delta_m(A_m(w)) = 0\}
 = 
\prod_{m \in \mathcal{M}_0(\Pr)}
\left[
1 - A_m(w)
\right],
\end{displaymath}
we simply replace the second $=$ sign by $\ge$ and then the proof of the theorem
goes through.

\subsection{Of Theorem \ref{theo-FDR control}}
\label{subsection-Proof of Second Theorem}

\begin{proof}
This proof closely mimics that of Theorem 7.1 in \cite{PenHabWu10} as presented in
\cite{PenHabWu10Supp}. As an aside, we mention that the seed of the {\em idea} of providing
a {\em class} of FDR-controlling multiple decision functions was planted upon
the realization that the proof of this Theorem 7.1 is independent of the choice of the
multiple decision size function.

The case with $q =0$ is trivial since then $\alpha^*(0) = 0$, so that
$F(\alpha^*(0)) = 0$. Thus we restrict to $q \in (0,1]$. By the defining property of
$\alpha^*(q)$ given in (\ref{alpha-star}),
we have that
\begin{equation}
\label{inequality at alpha*}
S(\alpha^*(q)) \ge \frac{1}{q} A_\bullet(\alpha^*(q))
\end{equation}
where $A_\bullet(\alpha) = \sum_{m \in \mathcal{M}} A_m(\alpha)$. Consequently, from
(\ref{F}),
\begin{equation}
\label{inequality for F(alpha*)}
F(\alpha^*(q))  \le q \frac{S_0(\alpha^*(q))}{A_\bullet(\alpha^*(q))} I\{S(\alpha^*(q) > 0\}
\le q \frac{S_0(\alpha^*(q))}{A_\bullet(\alpha^*(q))}.
\end{equation}

For $\alpha \in [0,1]$, define the sub-$\sigma$-field
\begin{equation}
\label{sub sigma field}
\mathcal{F}_\alpha \equiv
\mathcal{F}_\alpha(\mathbf{\Delta},\mathbf{A}) =
\sigma\left\{
\delta_m(A_m(\beta)):\ \beta \in [\alpha,1], m \in \mathcal{M}
\right\}.
\end{equation}
Observe that $\mathfrak{F} = (\mathcal{F}_\alpha: \ \alpha \in [0,1])$ is a
decreasing collection of sub-$\sigma$-fields of $\mathcal{F}$. By its definition
$\alpha^*(q)$ is an $\mathfrak{F}$-stopping time.

Let us define the process $\mathbf{T}_0 = (T_0(\alpha):\ \alpha \in [0,1])$ according to
\begin{displaymath}
T_0(\alpha) \equiv T_0(\alpha;\mathbf{\Delta},\mathbf{A},\Pr) =
\sum_{m \in \mathcal{M}_0(\Pr)}
\frac{\delta_m(A_m(\alpha))}{A_m(\alpha)}.
\end{displaymath}
Fix $0 \le \alpha \le \beta \le 1$. Then, since
$\delta_m \in \{0,1\}$, we have
\begin{eqnarray*}
\lefteqn{ E_\Pr\{T_0(\alpha) | \mathcal{F}_\beta\} } \\ & = &
\sum_{m \in \mathcal{M}_0(\Pr)}
E_\Pr\left\{\frac{\delta_m(A_m(\alpha))}{A_m(\alpha)} | \mathcal{F}_\beta\right\} \\
& = & \sum_{m \in \mathcal{M}_0(\Pr)} \left[\frac{1}{A_m(\alpha)}\right]
\Pr\{\delta_m(A_m(\beta)) = 1|\mathcal{F}_\beta\} \times \\ &&
E_\Pr\left\{{\delta_m(A_m(\alpha))} | \delta_m(A_m(\beta)) = 1\right\} \\
& = &  \sum_{m \in \mathcal{M}_0(\Pr)}
\delta_m(A_m(\beta)) \frac{1}{A_m(\alpha)} \frac{A_m(\alpha)}{A_m(\beta)}, \ \mbox{a.e.}\ [\Pr] \\
& = & T_0(\beta).
\end{eqnarray*}
The second equality follows from (D3), whereas
the second-to-last equality follows since
\begin{eqnarray*}
\lefteqn{ E_\Pr\left\{{\delta_m(A_m(\alpha))} 
| \delta_m(A_m(\beta)) = 1\right\} }
\\
& = &
\frac{\Pr\{\delta_m(A_m(\alpha))=1,\delta_m(A_m(\beta))=1\}}
{\Pr\{\delta_m(A_m(\beta)) = 1\}} \\
& = &
\frac{\Pr\{\delta_m(A_m(\alpha))=1\}}
{\Pr\{\delta_m(A_m(\beta)) = 1\}} \\
& = & \frac{A_m(\alpha)}{A_m(\beta)}
\end{eqnarray*}
because of condition (A2) for the $A_m(\cdot)$s and conditions (D2) and (D4) 
for the $\delta_m(\cdot)$s.
The above results show that, under $\Pr$,
\begin{displaymath}
\left\{\left(T_0(\alpha),\mathcal{F}_\alpha\right):\ \alpha \in [0,1]\right\}
\end{displaymath}
forms a reverse martingale process. Further, observe that $T_0(1) = |\mathcal{M}_0(\Pr)|$
a.e. $[\Pr]$ due to conditions (D1) and (A1). Thus, 
\begin{displaymath}
E_\Pr(T_0(1)) = |\mathcal{M}_0(\Pr)|.
\end{displaymath}

From the inequality in (\ref{inequality for F(alpha*)}), we obtain
\begin{eqnarray*}
E_\Pr[F(\alpha^*(q))] & \le &
q E_\Pr\left[\frac{S_0(\alpha^*(q))}{A_\bullet(\alpha^*(q))}\right] \\
& = & q \sum_{m \in \mathcal{M}_0(\Pr)} E_\Pr\left[
\frac{\delta_m(\alpha^*(q))}{A_\bullet(\alpha^*(q))} \right] \\
& = & q \sum_{m \in \mathcal{M}_0(\Pr)} E_\Pr\left[
\frac{\delta_m(\alpha^*(q))}{A_m(\alpha^*(q))}
\frac{A_m(\alpha^*(q))}{A_\bullet(\alpha^*(q))} \right] \\
& \le & q \left[ \sup_{\alpha \in [0,1]} \max_{m \in \mathcal{M}_0(\Pr)}
\frac{A_m(\alpha)}{A_\bullet(\alpha)} \right] E_\Pr\left[T_0(\alpha^*(q))\right] \\
& \le & q \frac{1}{|\mathcal{M}_0(\Pr)|} E_\Pr[T_0(1)] \\
& = & q \frac{|\mathcal{M}_0(\Pr)|}{|\mathcal{M}_0(\Pr)|} \\
& = & q,
\end{eqnarray*}
where the last inequality is obtained using condition (A4) and by invoking the
Optional Sampling Theorem for (reverse) martingales (cf., \cite{Doo53}),
and the second-to-last equality
because of $E_\Pr[T_0(1)] = |\mathcal{M}_0(\Pr)|$.

Note that, in particular, since the \Sidak\ multiple decision size function
$\mathbf{A}^S$ always satisfies condition (A4) for {\em all}
$\Pr \in \mathcal{P}$, then $\forall \mathbf{\Delta} \in \mathfrak{D},
\forall \Pr \in \mathcal{P}$, we have the property
\begin{equation}
\label{Sidak FDR property}
E_\Pr\left\{
F(\alpha^*(q;\mathbf{\Delta},\mathbf{A}^S); \mathbf{\Delta},\mathbf{A}^S)
\right\} \le q.
\end{equation}

Let us denote by $\mathcal{P}_0 = \{\Pr \in \mathcal{P}:\ \mathcal{M}_0(\Pr) =
\mathcal{M}\}$. Observe that for $\Pr \in \mathcal{P}_0$, condition (A4) will
not be satisfied unless the multiple decision size function
$\mathbf{A}$ has identical components,
in essence, a \Sidak\ multiple decision size
function form. We still therefore need to
establish that for an arbitrary
$\mathbf{A} \in \mathfrak{S}$ and a $\Pr \in \mathcal{P}_0$,
\begin{displaymath}
E_\Pr\left\{
F(\alpha^*(q;\mathbf{\Delta},\mathbf{A}); \mathbf{\Delta},\mathbf{A})
\right\} \le q.
\end{displaymath}
For such a $\Pr \in \mathcal{P}_0$, we have
$F(\alpha;\mathbf{\Delta},\mathbf{A}) =
I\{S(\alpha;\mathbf{\Delta},\mathbf{A}) > 0\}$, so that
\begin{eqnarray*}
\lefteqn{ E_\Pr[F(\alpha^*(q;\mathbf{\Delta},\mathbf{A});\mathbf{\Delta},\mathbf{A})] } \\ &  = &
\Pr\{S(\alpha^*(q;\mathbf{\Delta},\mathbf{A});\mathbf{\Delta},\mathbf{A}) > 0\} \\ & = &
\Pr\{\alpha^*(q;\mathbf{\Delta},\mathbf{A}) > 0\}.
\end{eqnarray*}
We have, for any $\mathbf{\Delta} \in \mathfrak{D}$ and any $\mathbf{A} \in
\mathfrak{S}$, that
\begin{equation}
\label{alpha* greater zero}
\left\{\alpha^*(q;\mathbf{\Delta},\mathbf{A}) > 0\right\} =
\bigcup_{\alpha \in (0,1]}
\left\{
\frac{S(\alpha; \mathbf{\Delta},\mathbf{A})}{A_\bullet(\alpha)} \ge \frac{1}{q}
\right\}.
\end{equation}

In Lemma D.1 of \cite{PenHabWu10Supp} it was established, using an inequality of
\cite{Hoe56}, that for $W_m(\eta_m), m \in \mathcal{M}$,
independent Bernoulli($\eta_m$) random variables with $\eta_m \in [0,1]$ and satisfying
$\prod_{m \in \mathcal{M}}(1-\eta_m) = 1-\alpha$, for each $t \ge 1$,
\begin{equation}
\label{Bernoulli inequality}
\Pr\left\{\frac{\sum_{m \in \mathcal{M}} W_m(\eta_m)}
{\sum_{m \in \mathcal{M}} \eta_m} \ge t \right\} \le
\Pr\left\{\frac{\sum_{m \in \mathcal{M}} W_m(\eta_m^S)}
{\sum_{m \in \mathcal{M}} \eta_m^S} \ge t \right\},
\end{equation}
where $\eta_m^S = 1-(1-\alpha)^{1/M}, m \in \mathcal{M}$.

Noting that, under $\Pr \in \mathcal{P}_0$, $\delta_m(A_m(\alpha))$s are independent
Bernoulli($A_m(\alpha)$), then by using the inequality in
(\ref{Bernoulli inequality}) and condition (A3), it follows
that for $q \in (0,1]$,
\begin{equation}
\label{prob inequality}
\Pr\left\{
\frac{S(\alpha;\mathbf{\Delta},\mathbf{A})}
{A_\bullet(\alpha)} \ge \frac{1}{q}\right\} \le
\Pr\left\{
\frac{S(\alpha;\mathbf{\Delta},\mathbf{A}^{S+})}
{A_\bullet^{S+}(\alpha)} \ge \frac{1}{q}\right\},
\end{equation}
where the \Sidak\ sizes $\mathbf{A}^{S+} = (A_m^{S+}, m \in \mathcal{M})$ in
(\ref{prob inequality}) have components
$$A_m^{S+} = 1 - (1 - \alpha^+)^{1/M}, m \in \mathcal{M},$$
with $\alpha^+$ satisfying
%
$\prod_{m \in \mathcal{M}} [1 - A_m(\alpha)] = 1 - \alpha^+.$
%
Observe that by (A3), we have necessarily that $\alpha^+ \le \alpha$.
Combining the results in (\ref{alpha* greater zero}) and (\ref{prob inequality}),
we obtain
\begin{displaymath}
\Pr\{\alpha^*(q;\mathbf{\Delta},\mathbf{A}) > 0\} \le
\Pr\{\alpha^*(q;\mathbf{\Delta},\mathbf{A}^{S+}) > 0\}.
\end{displaymath}
But since we have already established that, for $\Pr \in \mathcal{P}_0$, we have
$$\Pr\{\alpha^*(q;\mathbf{\Delta},\mathbf{A}^{S+}) > 0\} \le q,$$
then it follows that
$\Pr\{\alpha^*(q;\mathbf{\Delta},\mathbf{A}) > 0\} \le q.$
This implies finally that
$$E_\Pr\{F(\alpha^*(q;\mathbf{\Delta},\mathbf{A});\mathbf{\Delta},\mathbf{A})\}
\le q$$ 
for any $\Pr \in \mathcal{P}_0$, thereby completing the proof of 
Theorem \ref{theo-FDR control}.
\end{proof}

In contrast to Theorem \ref{theo-FWER control} where we were able to have the weaker
version of condition (D4) given in (\ref{weaker D4}), we could not do this
for Theorem \ref{theo-FDR control}. The reason is that we could {\em not} 
conclude under this weaker condition
that the process $\{(T_0(\alpha),\mathcal{F}_\alpha): \alpha \in [0,1]\}$ is a reverse
supermartingale, which would have allowed us to get the desired result. 
It may be possible that
under certain situations we do have this supermartingale property, but the weaker
condition (\ref{weaker D4}) appears not sufficient for this property to hold
in general.

\section{Representations of MDFs in Terms of the Generalized $P$-Values}
\label{section-Gen P Value}

This section expresses the MDFs $\delta^\dagger(q)$
in (\ref{FWER-controlling MDF}) and $\delta^*(q)$ in
(\ref{FDR-control MDF}) in terms of the
generalized $p$-value statistics defined in (\ref{generalized p-value}).
Define the anti-rank statistics vector via
\begin{equation}
\label{anti-ranks}
((1),(2),\ldots,(M)): (\mathcal{X},\mathcal{F}) \rightarrow (\mathfrak{M},\sigma(\mathfrak{M}))
\end{equation}
where $\mathfrak{M}$ is the space of all possible permutations of $\mathcal{M}$, and such
that
\begin{displaymath}
\alpha_{(1)} < \alpha_{(2)} < \ldots < \alpha_{(M)}.
\end{displaymath}

Let us first consider the random variable $\alpha^\dagger(q)$ in
(\ref{alpha-dag}). We see from its definition and those of the generalized
$p$-value statistics that, for some 
$J \in \bar{\mathcal{M}} \equiv \{0\} \cup \mathcal{M}$,
we have
\begin{displaymath}
\alpha^\dagger(q) \in [\alpha_{(J)}, \alpha_{(J+1)})
\end{displaymath}
if and only if
\begin{eqnarray*}
& \forall j \in \{1,2,\ldots,J\}: \prod_{m \in \mathcal{M}}
[1 - A_{(m)}(\alpha_{(j)}]^{1 - \delta_{(m)}[A_{(m)}(\alpha_{(j)})-]} \ge 1 - q; & \\
& \prod_{m \in \mathcal{M}}
[1 - A_{(m)}(\alpha_{(J+1)}]^{1 - \delta_{(m)}[A_{(m)}(\alpha_{(J+1)})-]} < 1 - q. &
\end{eqnarray*}
From the definition of the generalized $p$-value statistics we further have
\begin{displaymath}
\delta_{(m)}[A_{(m)}(\alpha_{(j)})-] = I\{m \le j-1\}.
\end{displaymath}
Consequently, by defining the $\bar{\mathcal{M}}$-valued random variable
\begin{equation}
\label{J dagger}
J^\dagger(q) =
\max
\left\{
k \in \mathcal{M}:\
\prod_{m=j}^M [1 - A_{(m)}(\alpha_{(j)})] \ge 1 - q, j=1,2,\ldots,k
\right\},
\end{equation}
we have the result that
\begin{displaymath}
\alpha^\dagger(q) \in \left[
\alpha_{(J^\dagger(q))}, \alpha_{(J^\dagger(q)+1)}
\right).
\end{displaymath}
As a consequence we obtain the representation of $\delta^\dagger(q)$
in (\ref{FWER-controlling MDF}) in terms of the $\alpha_m$s given by
\begin{equation}
\label{equiv delta dagger}
\delta^\dagger(q) \equiv
\left(
\delta_m(A_m(\alpha^\dagger(q))), m \in \mathcal{M}
\right)
=
\left(
\delta_m(A_m(\alpha_{(J^\dagger(q))})), m \in \mathcal{M}
\right),
\end{equation}
where we used the fact that, for each $m \in \mathcal{M}$,
$\delta_m$ is constant in the interval
\begin{displaymath}
[A_m(\alpha_{(J^\dagger(q))}), A_m(\alpha_{(J^\dagger(q)+1)}).
\end{displaymath}

Next let us consider the random variable $\alpha^*(q)$ in
(\ref{alpha-star}). We may re-express its defining equation via
\begin{displaymath}
\alpha^*(q) = \sup\left\{
\alpha \in [0,1]:
\sum_{m \in \mathcal{M}} A_{(m)}(\alpha) \le q \sum_{m \in \mathcal{M}}
\delta_{(m)}[A_{(m)}(\alpha)]
\right\}.
\end{displaymath}
But, since $\sum_{m \in \mathcal{M}} \delta_{(m)}[A_{(m)}(\alpha_{(j)})] = j$, then
\begin{displaymath}
\alpha^*(q) \in [\alpha_{(J)}, \alpha_{(J+1)})
\end{displaymath}
if and only if
\begin{eqnarray*}
& \sum_{m \in \mathcal{M}} A_{(m)}(\alpha_{(J)}) \le qJ; & \\
& \forall j \in \{J+1,J+2,\ldots,M\}: \sum_{m \in \mathcal{M}} A_{(m)}(\alpha_{(j)}) > qj. &
\end{eqnarray*}
Defining the $\bar{\mathcal{M}}$-valued random variable
\begin{equation}
\label{J star}
J^*(q) = {\max}\left\{k \in \mathcal{M}:\
\sum_{m \in \mathcal{M}} A_{(m)}(\alpha_{(k)}) \le qk \right\},
\end{equation}
we then have that 
\begin{displaymath}
\alpha^*(q) \in \left[\alpha_{(J^*(q))}, \alpha_{(J^*(q)+ 1)}\right).
\end{displaymath}
As a consequence, an equivalent representation of the MDF $\delta^*(q)$ in
(\ref{FDR-control MDF}) in terms of the $\alpha_m$s is provided by
\begin{equation}
\label{equiv delta star}
\delta^*(q) \equiv
\left(
\delta_m(A_m(\alpha^*(q))), m \in \mathcal{M}
\right)
=
\left(
\delta_m(A_m(\alpha_{(J^*(q))})), m \in \mathcal{M}
\right).
\end{equation}
The representations in (\ref{equiv delta dagger}) for $\delta^\dagger(q)$ and
(\ref{equiv delta star}) for $\delta^*(q)$ provide alternative computational 
approaches since, instead of computing
$\alpha^\dagger(q)$ and $\alpha^*(q)$, we may simply compute the 
generalized $p$-values, then $J^\dagger(q)$ and $J^*(q)$, and then finally
the realizations of the decision functions.

For a simple application, let us see what becomes of the MDFs $\delta^\dagger(q)$
and $\delta^*(q)$ if we use the \Sidak\ multiple decision size function $\mathbf{A}^S$
given in (\ref{Sidak size}). We use the alternate representations just obtained above.
By simple manipulations, we immediately obtain that
\begin{eqnarray*}
& J^\dagger(q)  = 
\max\left\{ k \in \mathcal{M}:\ 
\alpha_{(j)} \le 1 - (1 - q)^{M/(M-j+1)}, j=1,2,\ldots,k
\right\}; & \\
& J^*(q)  =  \max\left\{ k \in \mathcal{M}:\ 
M[1 - (1 - \alpha_{(k)})]^{1/M} \le qk
\right\}. &
\end{eqnarray*}
But, for these \Sidak\ size functions, the (ordinary) $p$-value statistics are
given by
\begin{displaymath}
P_m = A_m^S(\alpha_m) = 1 - (1 - \alpha_m)^{1/M}, m \in \mathcal{M}.
\end{displaymath}
Re-expressing the $J^\dagger(q)$ and $J^*(q)$ in terms of these $p$-values,
we easily obtain by simple manipulations that
\begin{eqnarray}
& J^\dagger(q)  = 
\max\left\{ k \in \mathcal{M}:\ 
P_{(j)} \le 1 - (1 - q)^{1/(M-j+1)}, j=1,2,\ldots,k
\right\};  \label{J dag Sidak} & \\
& J^*(q)  =  \max\left\{ k \in \mathcal{M}:\ 
P_{(k)} \le {qk}/{M}
\right\}. \label{J star Sidak} & 
\end{eqnarray}
Observe that $J^\dagger(q)$ in (\ref{J dag Sidak}) leads to the
step-down sequential \Sidak\ FWER-controlling procedure, see \cite{DudLaa08}; 
whereas, $J^*(q)$ in (\ref{J star Sidak})
is the usual form of the step-up Benjamini-Hochberg FDR-controlling 
procedure in \cite{BenHoc95}.
Thus, through the \Sidak\ sizes, we are able to obtain from our
formulation two popular MDFs for FWER and FDR control as special cases
of the MDFs $\delta^\dagger(q)$ and $\delta^*(q)$!

\section{Towards the Development of Optimal MDFs}
\label{section-Optimal MDFs}

Finally, in this subsection, we indicate, 
without going into much detail, the potential utility of the classes of MDFs
arising from Theorems \ref{theo-FWER control} and \ref{theo-FDR control} in the context of
obtaining MDFs with some optimality properties,
especially in non-exchangeable multiple hypotheses testing settings,
which are those where the power characteristics of the $M$ test functions
are not identical.
 
Let us fix a multiple decision process $\mathbf{\Delta}
\in \mathfrak{D}$ and fix a probability measure $\Pr_1 \in \mathcal{P}$. 
Define the mappings
$\pi_m: [0,1] \rightarrow [0,1]$ for $m \in \mathcal{M}$ according to
\begin{equation}
\label{ROCs}
\pi_m(\alpha;\Pr_1) = E_{\Pr_1} [\delta_m(\alpha)],\ \alpha \in [0,1].
\end{equation}
When viewed as a function of $\Pr_1$, $\pi_m(\alpha;\cdot)$ is the power function of $\delta_m$
when it is allocated a size of $\alpha$. Of interest to us, though, is to view it as a function
of $\alpha$ for the fixed $\Pr_1$. In this case, $\pi_m(\cdot;\Pr_1)$ is the receiver operating
characteristic curve (ROC) of the $m$th test function. Assume that for each $m \in \mathcal{M}$,
the mapping $\alpha \mapsto \pi_m(\alpha;\Pr_1)$
is strictly increasing with $\pi_m(1;\Pr_1) = 1$ and twice-differentiable.

Suppose it is desired to strongly control the overall FWER or FDR at some
pre-specified level $q \in [0,1]$, but at the same
time maximize the total (or average) power at $\Pr = \Pr_1$.
Our idea, partly implemented in
\cite{PenHabWu10}, is to first obtain the optimal multiple decision
size function for {\em weak} FWER control
associated with $\mathbf{\Delta}$, denoted by
$\mathbf{A}^* = (A_m^*(\alpha), m \in \mathcal{M}) \in \mathfrak{S}$. 
This is the multiple decision size function $\mathbf{A}$ satisfying the condition
\begin{displaymath}
\forall \alpha \in [0,1]:\ \prod_{m \in \mathcal{M}} [1 - A_m(\alpha)] = 1 - \alpha,
\end{displaymath}
and such that the total power 
at $\Pr = \Pr_1$, given by $\sum_{m \in \mathcal{M}} \pi_m(A_m(\alpha); \Pr_1)$, is
maximized. Under regularity conditions on the ROC functions,
the optimal $\mathbf{A}^*$ function could be
obtained using Lagrangian optimization, for instance, see
Theorem 4.3 in \cite{PenHabWu10} which is an implementation when the
individual test functions coincide with the Neyman-Pearson
most powerful tests.

Now, having determined the optimal multiple decision size function 
$\mathbf{A}^*$ associated with $\mathbf{\Delta}$,
which is at this point optimal only in
the sense of {\em weak} FWER control, we can then apply 
Theorem \ref{theo-FWER control} to obtain the MDF
$\delta^\dagger(q;\mathbf{\Delta},\mathbf{A}^*)$ 
which will {\em strongly} control the FWER at $q$; or apply
Theorem \ref{theo-FDR control} to obtain the MDF 
$\delta^*(q;\mathbf{\Delta},\mathbf{A}^*)$ which
will control the FDR at $q$. 

By virtue of the choice of the size process $\mathbf{A}^*$, which is
tied-in to the multiple decision process $\mathbf{\Delta}$ and the target 
probability measure $\Pr_1$, we expect
that the MDFs $\delta^\dagger(q;\mathbf{\Delta},\mathbf{A}^*)$ 
and $\delta^*(q;\mathbf{\Delta},\mathbf{A}^*)$
will perform better with respect to overall power at $\Pr_1$ relative to,
for example, the sequential
\Sidak\ MDF or the BH MDF, which we saw from the preceding section are MDFs
arising from the
\Sidak\ multiple decision size function, a size function that may not be
optimal for the chosen multiple decision size process $\mathbf{\Delta}$.
For instance, results of a modest simulation study in \cite{PenHabWu10}
demonstrated the improvement over the BH procedure of the MDF $\delta^*$ in a
specific setting.
Further improvements in power performances could
be achieved by proper choice of the multiple decision process $\mathbf{\Delta}$,
such as, for example, choosing it
to have components that are uniformly most powerful (UMP) or
uniformly most powerful unbiased (UMPU)
test functions. These issues, however, will be deferred for future work,
but we expect that the classes of MDFs presented here will play a
central role in dealing with these more complex 
multiple hypotheses testing problems. 

We close by pointing out that other approaches have also been proposed for
obtaining MDFs possessing certain optimality properties. 
Relevant papers pertaining to optimality
are \cite{WesKriYou98, WesKri01, 
GenWas03, GenWasRoe06, Sto07, SunCai07, SarZhoGho08, KanYeLiuAllGao09, RoqWie09}. 
Procedures with a Bayes or an empirical Bayes flavor
can be found in \cite{MulParRobRou04, ScoBer06, Efr08, GuiMulZha09}. 
In addition, it is also of interest to extend our
results to settings where the components of $\{\delta_m: m \in \mathcal{M}_0(\Pr)\}$ 
are dependent as in \cite{SarCha97, BenYek01}; 
see also the review article \cite{SarSANKHYA08}.
Another possible extension is to consider generalized FWER and FDR
as in \cite{Sar07}. However, we defer consideration of such extensions 
for future work.

\section*{Acknowledgements}

The authors are grateful to Professor Sanat Sarkar for many discusssions which were 
highly beneficial to this work.

\baselineskip=12pt

\bibliographystyle{acmtrans-ims}

\end{document}